\documentclass{amsart}
\usepackage{amsmath, amsthm, amsfonts}

\newtheorem{thm}{Theorem}[section]

\def\iprod{\otimes} 

\newcommand{\rsc}[2]{\left\langle#1, #2\right\rangle} 

\newcommand\bN{\mathbb N}

\newcommand\Q{\mathbb Q}

\newcommand\kp{*}

\def\tpa{{\sc itensor\_de}}

\def\bK{\mathbb K}
\def\bN{\mathbb N}
\def\bQ{\mathbb Q}

\def\cS{\mathcal S}

\newcommand{\mm}[1]{\mathsf{#1}}
\newcommand{\HH}{\mathsf{H}}
\newcommand{\Sh}{\mathsf{S}}
\newcommand{\EE}{\mathsf{E}}
\newcommand{\MM}{\mathsf{M}}

\def\tpa{{\sc itensor\_de}}

\author[Marni Mishna]{Marni Mishna\\ Dept. Mathematics, Simon Fraser University}
\address{{Dept. Mathematics,}\\ { Simon Fraser University,}\\
 { Burnaby, Canada}\\ { V5A 1S6}}
\email{mmishna@sfu.ca {\em Webpage:} http://math.sfu.ca/\~{}mmishna}

\title{Kronecker product identities from D-finite symmetric functions}
\thanks{This work was supported in part by NSERC}
\begin{document}
\begin{abstract}
  Using an algorithm for computing the symmetric function Kronecker
  product of D-finite symmetric functions we find some new Kronecker
  product identities. The identities give closed form formulas for
  trace-like values of the Kronecker product.
\end{abstract}

\maketitle

\section*{Introduction}
In the process of showing how the scalar product of symmetric
functions can be used for enumeration purposes,
Gessel~\cite{Gessel90}, proved that this product, and the Kronecker
product, {\em preserve D-finiteness}. Roughly, this means that if~$F$
and~$G$ are symmetric functions which both satisfy a particular kind
of system of linear differential equations, then so will the scalar
and Kronecker products of these functions. In an earlier work~\cite{ChMiSa05},
we give algorithms to calculate both of these systems of differential
equations.

In this short note we use this algorithm in a symbolic way to find explicit
expressions for Kronecker products of pairs of several common series of
symmetric functions, such as complete ($\HH=\sum_n h_n$), elementary
($\EE=\sum_n e_n$) and Schur ($\Sh=\sum_n\sum_{\lambda \vdash n} s_\lambda$).
Proposition 12 of ~\cite{ChMiSa05}, is the following identity,
\begin{equation*}
\left(\sum_\lambda s_\lambda\right)\otimes\left(\sum_\lambda
s_\lambda\right)={\exp\left(\sum_{n\geq
1}\frac{p_{2n-1}}{(2n-1)(1-p_{2n-1})}\right)}{\left({\prod_{n\geq
1}\left(1-p_n^2\right)}\right)^{-1/2}}.
\end{equation*}

This is the generating series of~$\sum_n\left(\sum_{\lambda\vdash
    n,\mu\vdash n, \lambda<\mu} s_\lambda\otimes s_\mu\right)$.  Here,
we apply the same technique to give a table of new identities of the
same flavour.

\section{Symmetric functions}
We use notation as in~Macdonald~\cite{Macdonald95} for our symmetric functions. 
A {\em partition} of a positive integer $n$ is
a decreasing sequence of integers $\lambda=(\lambda_1, \lambda_2,
\ldots, \lambda_k)$ whose sum is $n$. This is denoted $\lambda\vdash n$. A
partition is written in either vector or power notation, for
example $(7,7,4,4,1)=[1\,4^2\,7^2]$ are both partitions of 23.
A {\em symmetric function} is a sum of monomials in a some variable
set, that is invariant under any permutation of that variable
set. We can write any symmetric function as a sum of {\em monomial
symmetric functions}, defined for the variable set $\{x_1, x_2,
\ldots\}$ with respect to some partition $\lambda$ as
\begin{equation*}
m_\lambda:=\sum_{\sigma\in\cS_{\bN\setminus\{0\}}}
\left(r_1!\,r_2!\dotsm\right)^{-1}
x_{\sigma(1)}^{\lambda_1}\dotsm x_{\sigma(k)}^{\lambda_k}.
\end{equation*}
For example,
$m_{(3,2,2)}=x_1^3x_2^2x_3^2+x_3^3x_2^2x_1^2+x_4^3x_1^2x_3^2+\ldots$. We
also have the {\em elementary symmetric functions},
$e_n=m_{\langle1^n\rangle}$, and $e_\lambda=e_{\lambda_1}\dotsm
e_{\lambda_k}$; the {\em complete symmetric functions}
$h_n=\sum_{\lambda\vdash n} m_\lambda$, and
$h_\lambda=h_{\lambda_1}\dotsm h_{\lambda_k}$; and {\em power sum
  symmetric functions} $p_n=m_{(n)}=x_1^n+x_2^n+\ldots$,
$p_\lambda=p_{\lambda_1}\dotsm p_{\lambda_k}$. We postpone the
definition of the Schur symmetric functions to the next section, where
we shall be better equipped. Any of the
$h_\lambda$, $p_\lambda$,  $e_\lambda$, or $s_\lambda$ can form a $\Q$-basis of the
vector space $\Lambda$ of symmetric functions. We can also
view~$\Lambda$ as the ring~$\Q[p_1, p_2 \ldots]$ and, finally, we 
also work in the ring~$\hat\Lambda=\Q[[p_1, p_2, \ldots]]$.
 
\subsection{The scalar product of symmetric functions}
The ring of symmetric series is endowed with a scalar product defined
as a symmetric bilinear form such that the bases~$(h_\lambda)$
and~$(m_\lambda)$ are dual to each other:
\begin{equation}\label{eq:scalhm}
\rsc{m_\lambda}{h_\mu}=\delta_{\lambda\mu}.
\end{equation}
It turns out that
\[
\langle p_\lambda,p_\mu\rangle=z_\lambda\delta_{\lambda,\mu},
\]
with~$z_\lambda=(1^{r_1}r_1!)(2^{r_2}r_2!)\dotsm$ when
$\lambda=[1^{r_1}2^{r_2}\cdots]$.

The Schur basis is an orthonormal symmetric function basis under this
scalar product. In fact, Schur functions can be defined as the result
of applying the Gram-Schmidt process for orthogonalizing a basis,
applied to the monomial basis with the partitions ordered
lexicographically\footnote{In such an ordering, $1^n<1^{n-1}2<\cdots< n$.}.

\subsection{Plethysm of symmetric functions}
To conclude this brief recollection of symmetric functions, we
describe one type of composition that turns out to be quite useful
here: {\em plethysm}, written $f[g]$. We can most easily define it
using the power sum symmetric functions. It is defined by
$p_n[\psi(p_1,p_2,\dots)]=\psi(p_n,p_{2n},\dots)$, along with
$(\phi_1+c\phi_2)[\psi]=\phi_1[\psi]+c\phi_2[\psi]$ and
$(\phi_1\cdot\phi_2)[\psi]=\phi_1[\psi]\cdot\phi_2[\psi]$.

\subsection{The Kronecker product of symmetric functions}
In the ring of symmetric functions the usual polynomial multiplication
serves as a product, but there is also a second product which arises
from the connection between symmetric function and the characters of
the symmetric group. This product has several names, including the
{\em Kronecker product\/}, the tensor product and the internal
product. Although we mostly follow the notation of
Macdonald~\cite{Macdonald95} for most matters relating to symmetric
functions, we shall refer to it here as the Kronecker product, and
denote it by~$\kp$. It was first described by Redfield as the {\em cap
  product\/} of symmetric functions and was rediscovered by
Littlewood~\cite{Littlewood56}. This product can be defined in
representation theory pointwise product of characters, which
corresponds to tensor products of representations, however here we use
the following relation to the power sum symmetric functions, and
extend linearly:
\begin{equation}\label{eq:iprod}
p_\lambda\iprod p_\mu =\delta_{\lambda\mu}z_\lambda p_\lambda.
\end{equation} 

Calculating the connection coefficients
$\gamma_{\lambda,\mu}^{(\rho)}$ for the Kronecker product in the Schur
basis
\[ s_\lambda\iprod s_\mu = \sum_\rho \gamma_{\lambda,\mu}^{(\rho)}
s_\rho \] is also challenging, and quite interesting. There are some
combinatorial interpretations of $\gamma_{\lambda,\mu}^{(\rho)}$ which have obtained results when
$\lambda,\rho$ and $\mu$ are of a particular form, such as work of
Goupil and Schaeffer~\cite{GoSh98}, Rosas~\cite{Rosas01a}, or Chauve
and Goupil~\cite{ChGo05}. The interest originates from the
correspondence with irreducible representations,
\newcommand{\cV}{\mathcal{V}}
\begin{equation*}
\chi(\cV_\lambda \kp \cV_\mu) = \sum_\rho \gamma_{\lambda,\mu}^{(\rho)}\,
\chi(\cV_\rho).
\end{equation*} 
When $\lambda$, $\mu$ and $\rho$ are all partitions of $n$,
$\gamma_{\lambda,\mu}^{(\rho)}$ is the multiplicity of a character in
the representation. For more details, the reader is pointed towards
the text of Sagan~\cite{Sagan01}.

To compute $\gamma_{\lambda,\mu}^{(\rho)}$ using computer algebra
systems, one typically expands the symmetric function into the power
sum basis and then applies
\eqref{eq:iprod} to a pairwise comparison of terms. (For example, in
the SF package of Stembridge.) As we mentioned in the introduction, we
introduced a generating function approach~\cite{ChMiSa05}. The
algorithm in~\cite{ChMiSa05} that we use is called~\tpa, and a Maple
implementation on the author's web page is available.

A second approach, summarized in~\cite{ScThWy93}, uses a reduced
notation that allows calculations with series of the form $\sum_n
s_{(n,\lambda_2, \dots, \lambda_k)}z^n$, for fixed $\lambda_2, \dots,
\lambda_k$. These computations are quite efficient; far more so than
expanding the power-sum basis. 

\section{Applications of D-finite symmetric series}
The algorithms we use do not compute the products
directly, rather they determine differential equations satisfied by
the resulting function. The existence of such differential equations is a
consequence of the D-finite closure properties of the scalar
product.  A function~$\phi$ is said to be {\em D-finite}
in $\bK[[x_1,\dots,x_r]]$ if and only if
the partial derivatives
$\partial_1^{\alpha_1}\dotsm\partial_r^{\alpha_r}\phi$ generate a
finite dimensional vector space over~$\bK(x_1,\dots,x_r)$. In this
case, $\phi$~is determined by a system of linear differential
equations.
  
In order to treat symmetric functions, however, we must consider
functions with an infinite number of variables. The function
$\phi(x_1, x_2, \ldots)$ is D-finite in $K[[x_1, x_2, \ldots]]$ if for
all~$r$, $\phi(x_1, \ldots, x_r, 0, \cdots)$ is D-finite in $K[[x_1,
x_2, \ldots, x_r]]$.  This case does not enjoy all of the closure
properties of the previous, nonetheless we have closure under $+$,
$\times$, $\partial_i$, extension of coefficients, rational
substitution, and exponentials of polynomials. We say that a symmetric
function $\phi\in\hat\Lambda$ is D-finite if it is D-finite in
$\bQ[[p_1,p_2,\dots]]$. For example, under this definition the two
following famous symmetric function sums $\HH$ and $\EE$  which we
introduced earlier satisfy the following relations,
\[
\HH=\exp\left(\sum_n\frac{p_n}{n}\right)
\text{ and }
\EE=\exp\left(\sum_n(-1)^n\frac{p_n}{n}\right),
\]
and thus are both D-finite.

It was Gessel~\cite{Gessel90} that first showed that the scalar product and the
Kronecker product both preserve D-finiteness. The work~\cite{ChMiSa05}
makes this effective by transforming the system
of differential equations satisfied by $F$ and $G$ in to one satisfied
by $F\iprod G$, or $\langle F,G\rangle$.

\subsection{Kronecker product calculations}
Many interesting problems which use the Kronecker product involve
symmetric functions, which once they are expressed in the power sum
basis, require an infinite number of~$p_n$. Thus, at first glance they
are seemingly unsuitable for direct application of our algorithms
which, after all, require finite input! One approach is to apply these
algorithms for several truncations of the symmetric functions and
generate information upon which reasonable conjectures can be
formulated. For each of these, we render the problem applicable by
setting most $p_n$'s to 0. That is, we solve a sequence of problems
involving an increasing number of $p_n$, and hope to identify a
pattern.

However, far more satisfying are the cases where there is sufficient
form and structure which can be exploited to find exact results. We
shall be more specific about precisely the ``form and structure'' we
can exploit in a moment. First we remark that one important such class
comes from symmetric series arising from plethysms. In this case, we
can reduce the Kronecker product of functions each with an infinite
number of $p_n$ variables to a finite number of symbolic calculations.

For example, if two symmetric functions~$F$ and~$G$ can be expressed
respectively  in the form 
\[F(p_1, p_2, \ldots)=\prod_{n\ge 1} f_n(p_n) \qquad \text{ and } \qquad
G(p_1, p_2, \ldots)=\prod_{n\ge 1} g_n(p_n),\]
then one can
easily deduce that 
\begin{equation}\label{eqn:prodform}
F\kp G=\prod_{n\ge 1}  f_n(p_n)\kp g_n(p_n).
\end{equation}
Essentially this follows from the fact that the Kronecker product of
two power sum symmetric functions of differing order is 0.
If, furthermore, the $f_n$ and $g_n$ are such that one can describe
them in a finite way using D-finite functions, we can apply this
method. 

Series which arise as plethysms of the form~$\HH[u]$ or~$\EE[u]$,
where~$u$ is a polynomial in the $p_i$, are precisely of this
form. For example the sum of all Schur functions is of this type:
\[
\Sh=\sum_\lambda
s_\lambda=\HH[p_1+\frac{1}{2}p_1^{2}-\frac{1}{2}p_2]
         =\exp\left(\sum_n\frac{p_n^2}{2n}+\frac{p_{2n-1}}{2n-1}\right).
\]
Thus, 
\[
\Sh=\left(\prod_{n\text { even}} \exp\left(\frac{p_n^2}{2n}\right)\right)
    \left(\prod_{n\text{ odd}}  \exp\left(\frac{p_n^2}{2n}+\frac{p_n}{n}\right)\right)
\]
We assign~$f_n$ as follows
\[f_{2n}=\exp\left(\frac{p_{2n}^2}{4n}\right)\quad \text{ and }\quad
f_{2n-1}=\exp\left(\frac{p_{2n-1}^2}{2}+\frac{p_{2n-1}}{2n-1}\right).
\]
Thus, to compute $\Sh\iprod \Sh$, we compute in turn $g_{2n}=f_{2n}\kp
f_{2n}$, and $g_{2n+1}=f_{2n+1}\kp f_{2n+1}$. We find $g_{2n}$ by
determining the differential equation that it satisfies, using \tpa
adapted to handle a formal parameter. The adaptation amounts to
performing a scalar product with adjunction formula
$p^\diamond=n\partial$ for a formal parameter $n$. This gives
\begin{equation*}
(1-p_n^2)\frac{\partial g_n(p_n)}{\partial p_n}+p_ng_n(p_n)=0,
\qquad\text{for even~$n$}.
\end{equation*}
We then solve for $g_n$. We do likewise for odd $n$, and then the
identity in the introduction follows. 

Carbonara {\em et al.}~\cite{CaCaRe03} are interested in the trace of
$\sum_n\left(\sum_{\lambda\vdash n,\mu\vdash n, \lambda<\mu}
  s_\lambda\kp s_\mu\right)$ given by
$\sum_n\left(\sum_{\lambda\vdash n} s_\lambda\kp
  s_\lambda\right)$. It is not immediately clear to me if our method
could be adapted directly to this kind of calculation.

\subsection{A family of identities}
We now apply the above approach to create a number of different identities. 
The following table summarizes results. These formulas for $\HH, \EE, \Sh,
\Sh\EE^{-1}$ and $\Sh\HH^{-1}$ are all derived in
Macdonald~\cite{Macdonald95}:
\[
\begin{array}{ll}
\HH=\sum h_nt^n=\exp\left(\sum_n \frac{p_n}{n}\right)\\
\EE=\sum e_nt^n=\exp\left(\sum_n (-1)^{n+1}\frac{p_n}{n}\right)\\
\Sh=\sum_\lambda s_\lambda t^{|\lambda|}=\exp\left(\sum_n
\frac{p^2_nt^{2n}}{2n}+\frac{p_{2n-1}t^{2n-1}}{2n-1}\right)\\
\Sh\EE^{-1}=\sum_{\lambda\text{ all parts odd}} s_\lambda t^{|\lambda|}\\
\Sh\HH^{-1}=\sum_{\lambda'\text{all parts even}} s_\lambda
t^{|\lambda|}
\end{array}
\]

\begin{thm}\label{thm:table}
Given the above definitions for $\HH, \EE$ and $\Sh$. Then, there is 
the following multiplication table for the Kronecker product, 
\begin{center}\large
\framebox{\begin{tabular}{c|llllll}
$\kp$&    $\HH$&    $\EE$&    $\Sh$&    $\Sh\HH^{-1}$&      $\Sh\EE^{-1}$\\ \hline
$\HH$&    $\HH$&    $\EE$&    $\Sh$&    $\Sh\HH^{-1}$&      $\Sh\EE^{-1}$\\
$\EE$&    &$\HH$&    $\Sh$&    $\Sh\EE^{-1}$&      $\Sh\HH^{-1}$\\
$\Sh$&    &       &       $\mm G\MM_{odd}$& $\mm G\mm N$&       $\mm G\mm N$      \\
$\Sh\HH^{-1}$&&     &       &       $\mm G\MM_{even}$&    $\mm G\mm P$\\
$\Sh\EE^{-1}$&&     &       &       &               $\mm G\MM_{even}$\\
\end{tabular}
}.
\end{center}
The products are expressed in terms of the following: 
\[
\begin{array}{ll}
\MM_{odd (even)} = \exp\left(\sum_{n\;odd
(even)}\frac{p_{n}}{n(1-p_n)}\right) \\
\mm{N} =\exp\left(\sum_{n} \frac{p_{n}^2}{2n(1-p_n^2)}\right)\\
\mm{P} = \exp\left(\sum_{n\;even} \frac{p_{n}}{n(1+p_n)}\right)\\
\mm{G} = \prod_{n\geq 1}\left(1-p_n^2\right)^{-1/2}.
\end{array}
\]

\end{thm}
A Maple worksheet with the calculations behind the above table is
available at the author's website. We welcome all suggestions for
other series of interest. It would equally easy to treat plethysms of
the form $\HH[\phi]$ for some symmetric polynomial $\phi$. A preliminary
review of the work of Scharf, Thibon, and Wybourne, for
example~\cite{ScTh91, ScThWy93} suggests
that there may be more to do with series of the form~ $\sum_n
s_{(n,\lambda_2, \dots, \lambda_k)}z^n$, for fixed $\lambda_2, \dots,
\lambda_k$ if they can be shown to be D-finite.

We are able to compute, with this method, expressions satified by some
powers of $\mm{S}=\sum_\lambda s_\lambda$,  with respect to the
Kronecker product, for example $\mm{S}\kp\mm{S}\kp\mm{S}$, but these
result in differential equations which we are presently unable to
solve into explicit expressions. 

\section*{Conclusion}
The symbolic application of tensor product calculation yields, rather
easily, families of Kronecker product identities. It is possible that
these identities could be exploited for group theoretic gain, however,
this remains to be investigated, as does finding connections between
our formulas, and that of the trace co-characters.

\subsection*{Acknowledgments}
This work was initiated during a visit to Project Algorithms, in part
from discussions with Fr\'ed\'eric Chyzak, and was funded in part by
the NSERC (Canada).  Thanks are due also to Rosa Orellana for an
interesting discussion on the trace co-characters and an anonymous
referee that suggested several interesting references.


\end{document}